\documentclass[11pt, twoside, reqno]{amsart}
\usepackage{amsmath,amsthm}
\usepackage{amssymb,latexsym}
\usepackage{enumerate}
\usepackage{amsfonts}
\usepackage{mathrsfs}
\usepackage{amsmath, amsthm}
\usepackage{amssymb}
\usepackage{fancyhdr}
\usepackage{color}
\usepackage[colorlinks,
linkcolor=red,
anchorcolor=blue,
citecolor=green]{hyperref}
\newtheorem{theorem}{\bf Theorem}[section]
\newtheorem{lemma}[theorem]{\bf Lemma}
\newtheorem{prop}[theorem]{\bf Proposition}
\newtheorem{coro}[theorem]{\bf Corollary}

\numberwithin{equation}{section}

\begin{document}
	
	\baselineskip=17pt
	
	\renewcommand{\thefootnote}{\fnsymbol {footnote}}
	
	\title[weak type estimates for sparse operators]{Sharp off-diagonal weighted weak type estimates for sparse operators}
	
	\author[Qianjun He and Dunyan Yan]{Qianjun He
		\quad Dunyan Yan}
	
	\address{School of Mathematics, Graduate University, Chinese Academy of Sciences, Beijing 100049, China}
	\email{heqianjun16@mails.ucas.ac.cn}
	\email{ydunyan@ucas.ac.cn}
	\thanks{This work was supported by National Natural Science Foundation of China (Grant Nos. 11471309 and 11561062)}
	
	\begin{abstract}
		We prove sharp weak type weighted estimates for a class of sparse operators that includes majorants of standard singular integrals, fractional integral operators, and square functions. These bounds are knows to be sharp in many cases, and our main new result is the optimal bound $$[w]_{A_{p,q}}^{\frac{1}{q}}[w^{q}]_{A_{\infty}}^{\frac{1}{2}-\frac{1}{p}}\lesssim[w]_{A_{p,q}}^{\frac{1}{2}-\frac{\alpha}{d}}$$ for proper conditions which satisfy that three index $p$, $q$ and $\alpha$ ensure weak type norm of fractional square functions on $L^{q}(w^{q})$ with $p>2$.
	\end{abstract}

	\maketitle

	\vspace*{0.5mm}
	\section{Introduction}
	\setcounter{equation}{0}
	We study weighted inequalities for sparse operators, which can be defined by
	\begin{equation}\label{define-1}
	\mathcal{A}_{\alpha,\nu}^{\mathcal{S}}(f):=\bigg(\sum_{Q\in\mathcal{S}}\langle f\rangle_{\alpha,Q}^{\nu}{\bf{1}}_{Q}\bigg)^{\frac{1}{\nu}},\quad \langle f\rangle_{\alpha,Q}=\frac{1}{|Q|^{1-\frac{\alpha}{d}}}\int_{Q}f,
	\end{equation}
	where $\nu>0$, $0\leq\alpha<d$ and $\mathcal{S}$ is a sparse collection of dyadic cubes, i.e. all (dyadic) cubes $Q\in\mathcal{S}$, there exists $E_{Q}\subset Q$ which are pairwise disjoint and $|E_{Q}|\geq\gamma|Q|$ with $0<\gamma<1$. Note that $\langle f\rangle_{Q}$ denote $\langle f\rangle_{\alpha,Q}$ with $\alpha=0$. And so far it it know that the operator $\mathcal{A}_{\alpha,\nu}^{\mathcal{S}}$ dominate large classes of classical operators $T$, relying upon the sparse domination formula
	\begin{equation}\label{ sparse domination formula}
	|Tf(x)|\lesssim\sum_{i=1}^{N}\mathcal{A}_{\alpha,\nu}^{\mathcal{S}_{i}}(|f|)(x),
	\end{equation}
	where the collections $\mathcal{S}_{i}$ depend on the function $f$. For $\nu=1$ and $\nu=2$ with $\alpha=0$, $T$ becomes the Calder\'{o}n-Zygmund singular integrals \cite{La1,L3} and Littlewood-Paley square functions \cite{L2,L1}, respectively. Thus, the various norm inequalities that we prove for $\mathcal{A}_{\alpha,\nu}^{\mathcal{S}}$ immediately translate to corresponding estimates for these classses of classical operators.
	
	  A weight $w$ on $\mathbb{R}^{d}$ is a locally integrable function $w$: $\mathbb{R}^{d}\rightarrow (0,+\infty)$. The class of all $A_{\infty}$ weights consists of all weights $w$ for which their $A_{\infty}$ characteristic
	  $$[w]_{A_{\infty}}:=\sup_{Q}\frac{1}{w(Q)}\int_{Q}M({\bf{1}}_{Q}w)<\infty,$$
	  where $M$ is the Hardy-Littlewood maximal function and the suprema take over cubes of sides are parallel to the coordinate axes.
	
	  More precisely, we are concerned with quantifying the dependence of various weighted operator norms on a mixture of the  two weight $A_{p,q}^{\alpha}$ characteristic
	  $$[w,\sigma]_{A_{p,q}^{\alpha}}:=\sup_{Q\in\mathcal{S}}|Q|^{q(\frac{\alpha}{d}-1)}w(Q)\sigma(Q)^{\frac{q}{p^{\prime}}}<\infty.$$
	 The study of such maixed bounds was initiated in \cite{HP}. All our estimates will be stated in a dual-weight formulation, in which the classical one-weight off-diagonal case $A_{p,q}$ as defined below.
	
	 Since we dealing with sparse operators, we also consider the sparse versions of the weight characteristics, where the supremums above  are over dyadic cubes only. This is a standing convention throughout this paper without further notice.
	
	 Throughout this paper, $1<p,p^{\prime},q<\infty$, $p$ and $p^{\prime}$ are conjugate indices, i.e. $1/p+1/p^{\prime}=1$. Formally, we will also define $p=1$ as conjugate to $p^{\prime}=\infty$ and vice versa.
	
	Now, we formulate our main results as follows.
	
	 \begin{theorem}\label{weak estimate for sparse operator}
	 		Let $0<\nu<\infty$, $0\leq\alpha<d$ and $1<p\leq q<\infty$. Let $w,\sigma$ be a pair of weights. Then
	 		\begin{equation}\label{weak type bound for sparse operator}
	 			\|\mathcal{A}_{\alpha,\nu}^{\mathcal{S}}(\cdot\sigma)\|_{L^{p}(\sigma)\rightarrow L^{q,\infty}(w)}\lesssim [w,\sigma]_{A_{p,q}^{\alpha}}^{\frac{1}{q}}
	 			\left\{
	 			\begin{aligned}
	 			&[w]_{A_{\infty}}^{\frac{1}{\nu}(1-(\frac{\nu}{p})^{2})}[\sigma]_{A_{\infty}}^{\frac{1}{\nu}(\frac{\nu}{p})^{2}}, & \quad & p=q>\nu\,\,\,\,\text{and} \,\,\,\,\alpha>0,  \\
	 			&[\sigma]_{A_{\infty}}^{\frac{1}{q}}, & \quad & p\leq\nu\leq q,\\
	 			&[w]_{A_{\infty}}^{(\frac{1}{\nu}-\frac{1}{p})_{+}},  & \quad & \text{other case}.
	 			\end{aligned}
	 			\right.
	 		\end{equation}
	 		{\noindent}where $x_{+}:=\max(x,0)$ in the exponent. Here and below, we simplify case analysis by interpreting $[w]_{A_{\infty}}^{0}=1$, whether or not $[w]_{A_{\infty}}$ is finite.
	 \end{theorem}
	  
	 Lacey and Scurry  \cite{LS} provided a method to proof of the case $ q<\nu$ of Theorem $\ref{weak estimate for sparse operator}$, and we merely repeat their one-weight proof in the two-weight off-diagonal case. For $p>\nu$, we bound
	\begin{equation}\label{weak type sharp bounds}
	[w,\sigma]_{A_{p,q}^{\alpha}}^{\frac{1}{q}}[w^{q}]_{A_{\infty}}^{\frac{1}{\nu}-\frac{1}{p}}\lesssim[w]_{A_{p,q}}^{\frac{1}{q}}[w]_{A_{p,q}}^{\frac{1}{\nu}-\frac{1}{p}}=[w]_{A_{p,q}}^{\frac{1}{\nu}-\frac{\alpha}{d}}
	\end{equation}
	 is new even in the one weight case for $\frac{1}{q}+\frac{\alpha}{d}=\frac{1}{p}$. For $\nu\leq q\leq\frac{\nu}{1-\frac{\nu\alpha}{d}}$, we also obtain the bounds $[w]_{A_{p,q}}^{\frac{1}{q}}$  and it has an additional logarithmic factor, taking the form $(1+\log[w^{q}]_{A_{\infty}})^{\frac{1}{\nu}}$. This form bouds which will be proved in Section 4.
	
	 Theorem $\ref{weak estimate for sparse operator}$ include several known cases, the Sobolev type case $\frac{1}{q}+\frac{\alpha}{d}=\frac{1}{p}$ of these results, together with strong type estimate and multilinear extensions, can also be recovered from Fackler and Hyt\"{o}nen \cite{FH},  Zorin-Kranich \cite{Z1} the recent general framework, respectively.
	
	 For $\nu=1$ and $\alpha=0$, $\eqref{ sparse domination formula}$ holds for all Calder\'{o}n-Zygmund operators. Lerner \cite{L3} first prove the result, and  Lacey \cite{La1} give the most general version, with a simplified proof in the paper \cite{L4}. The bound $\eqref{weak type bound for sparse operator}$ in this case was obtained in \cite{HP} for $p=q=1$. In \cite{HL}, H\"{a}nninen and Lorist consider the sparse domination for the lattice Hardy-Littlewood maximal operator, and their obtained sharp weighted weak $L^{p}$ estimates.
	 
	  For $\nu=2$ and $\alpha=0$, $\eqref{ sparse domination formula}$ holds for several square function operators of Littlewood-Paley type \cite{DLR,LS,L2}. For $p = q$, the mixed bound $\eqref{weak type bound for sparse operator}$, even for general $\nu>0$, is from \cite{HL,LL}. This improves the pure $A_{p}$ bound of \cite{DLR,LS,L2}.
	  
	  	For $\nu=1$ and $0<\alpha<d$, $\eqref{ sparse domination formula}$ holds for the fractional integral operator \cite{LMPT}
	  	\begin{equation}\label{define fractional integral operator}
	  	I_{\alpha}f(x):=\int_{\mathbb{R}^{d}}\frac{f(y)}{|x-y|^{n-\alpha}}dy.
	  	\end{equation}
	  In the case for $p<q$,  $\eqref{weak type bound for sparse operator}$ are due to {\cite{CM1}}. The Sobolev type case with $\frac{1}{q}+\frac{\alpha}{d}=\frac{1}{p}$ was obtained
	  	by the same authors in \cite{CM}. Additional complications with $p=q$, which lead to the weaker version of our bound $\eqref{weak type bound for sparse operator}$, have been observed and addressed in different
	  	ways in \cite{CM1,CM}.
	
	For $\nu>0$ and $\alpha=0$, the bound $\eqref{weak type bound for sparse operator}$ in the case was obtained by Hyt\"{o}nen and Li \cite{HL} for $p=q\in (1,\infty)$.
	
	Theorem $\ref{weak estimate for sparse operator}$ with $\nu=2$ completes the picure of sharp weighted inequalities for fractional square functions, aside from the remaining case of $2\leq q\leq \frac{2}{1-\frac{2\alpha}{d}}$. Namely, $[w]_{A_{p,q}}^{\max(\frac{1}{q},\frac{1}{2}-\frac{\alpha}{d})}$ is the optimal bound among all possible bounds of form $\Phi([w]_{A_{p,q}})$ with an incrasing function $\Phi$. This was shown by Hyt\"{o}nen and Li \cite{HL},  Lacey and Scurry \cite{LS} in the category of power type function $\Phi(t)=t^{\beta}$; a variant of their argument proves the general claim, as we show in the last section.
	
	To prove the above results, we need the following characterization, which is essentially due to Lai \cite{Lai}; we supply the necessary details to cover the cases that were not explicitly treated in \cite{Lai}.
	
	\begin{theorem}\label{fractional sparse operator controlled by testing constants}
		Let $1<p\leq q<\infty$, $\nu>0$, $p>\nu$ and $0\leq\alpha<d$. Let $w$, $\sigma$ be a pair of weights. Then
		$$
			\|\mathcal{A}_{\alpha,\nu}^{\mathcal{S}}(\cdot\sigma)\|_{L^{p}(\sigma)\rightarrow L^{q,\infty}(w)}^{\nu}\simeq \mathscr{T}^{*},
		$$
		where the testing constants defined by
		$$
		\mathscr{T}^{*}:=\sup_{R\in\mathcal{S}}w(R)^{-\frac{1}{(\frac{q}{\nu})^{\prime}}}\big\|\mathop{\sum_{Q\in\mathcal{S}}}_{Q\subset R}\langle\sigma\rangle_{\alpha,Q}^{\nu-1}\langle w\rangle_{\alpha,Q}{\bf{1}}_{Q}\big\|_{L^{(\frac{p}{\nu})^{\prime}}(\sigma)}.
		$$
	\end{theorem}
	
The case $p>\nu$ of Theorem $\ref{weak estimate for sparse operator}$ is a consequence of Theorem $\ref{fractional sparse operator controlled by testing constants}$.  The estimation of the testing $\mathscr{T}^{*}$ given by Fackler and Hyt\"{o}nen \cite{FH} and their obtained following result.	
	
\begin{prop}\label{estiamte for two testing constants}
	Let $\nu>0$, $0\leq\alpha<d$, $p>\nu$ and $1<p\leq q<\infty$.  For $\mathscr{T}^{*}$ as in Theorem $\ref{fractional sparse operator controlled by testing constants}$, we have
	$$
	\mathscr{T}^{*}\lesssim
	[w,\sigma]_{A_{p,q}^{\alpha}}^{\frac{\nu}{q}}
	\left\{
	\begin{aligned}
	&[w]_{A_{\infty}}^{1-(\frac{\nu}{p})^{2}}[\sigma]_{A_{\infty}}^{(\frac{\nu}{p})^{2}}, & \quad & p=q\,\,\,\,\text{and} \,\,\,\,\alpha>0,  \\
	&[w]_{A_{\infty}}^{1-\frac{\nu}{p}},  & \quad & \text{other case}.
	\end{aligned}
	\right.
	$$
\end{prop}

The plan of the paper is as follows: We come with the proof of Theorem $\ref{fractional sparse operator controlled by testing constants}$, this completes the proof of Theorem $\ref{weak estimate for sparse operator}$ in the case of $p>\nu$. The remaining case of Theorem $\ref{weak estimate for sparse operator}$ for $p\leq \nu$ is then handled in Section 3. In the final scetion, we discuss the sharpness of our weak type estimates by modifying the example given by Lacey and Scurry \cite{LS}.

	 \section{Proof of Theorem $\ref{fractional sparse operator controlled by testing constants}$}

As mentioned, Theorem $\ref{fractional sparse operator controlled by testing constants}$ is essentially duo to Hyt\"{o}nen and Li \cite{HL}.

First, we give the following lemma.
\begin{lemma}\label{Lemma-1.1}
Let $w,\sigma$ be a pair of weights and $p>\nu>0$.
$$
\|	\mathcal{A}_{\alpha,\nu}^{\mathcal{S}}(\cdot \sigma)\|_{L^{p}(\sigma)\rightarrow L^{q,\infty}(w)}\simeq\sup_{\|f\|_{L^{p}(\sigma)}=1}\big\|\sum_{Q\in\mathcal{S}}\langle\sigma\rangle_{\alpha,Q}^{\nu}\langle f^{\nu}\rangle_{Q}^{\sigma}{\bf{1}}_{Q}\big\|_{L^{\frac{q}{\nu},\infty}(w)}
$$
\end{lemma}
$Proof$. By the definition of $\mathcal{A}_{\alpha,\nu}^{\mathcal{S}}$, we have
\begin{align*}
&\|\mathcal{A}_{\alpha,\nu}^{\mathcal{S}}(\cdot \sigma)\|_{L^{p}(\sigma)\rightarrow L^{q,\infty}(w)}=\sup_{\|f\|_{L^{p}(\sigma)}=1}\big\|\sum_{Q\in\mathcal{S}}\langle f\sigma\rangle_{\alpha,Q}^{\nu}{\bf{1}}_{Q}\big\|_{L^{\frac{q}{\nu},\infty}(w)}\\
&=\sup_{\|f\|_{L^{p}(\sigma)}=1}\big\|\sum_{Q\in\mathcal{S}}\langle \sigma\rangle_{\alpha,Q}^{\nu}(\langle f\rangle_{Q}^{\sigma})^{\nu}{\bf{1}}_{Q}\big\|_{L^{\frac{q}{\nu},\infty}(w)}\\
&\leq \sup_{\|f\|_{L^{p}(\sigma)}=1}\big\|\sum_{Q\in\mathcal{S}}\langle \sigma\rangle_{\alpha,Q}^{\nu}\langle (M_{\sigma}(f))^{\nu}\rangle_{Q}^{\sigma}{\bf{1}}_{Q}\big\|_{L^{\frac{q}{\nu},\infty}(w)}\\
&=\sup_{\|f\|_{L^{p}(\sigma)}=1}\big\|\sum_{Q\in\mathcal{S}}\langle \sigma\rangle_{\alpha,Q}^{\nu}\left\langle \left(\frac{M_{\sigma}(f)}{\|M_{\sigma}(f)\|_{L^{p}(\sigma)}}\right)^{\nu}\right\rangle_{Q}^{\sigma}{\bf{1}}_{Q}\big\|_{L^{\frac{q}{\nu},\infty}(w)}\|M_{\sigma}(f)\|_{L^{p}(\sigma)}^{\nu}\\
&\lesssim\sup_{\|g\|_{L^{p}(\sigma)}=1}\big\|\sum_{Q\in\mathcal{S}}\langle\sigma\rangle_{\alpha,Q}^{\nu}\langle g^{\nu}\rangle_{Q}^{\sigma}{\bf{1}}_{Q}\big\|_{L^{\frac{q}{\nu},\infty}(w)},
\end{align*}
where in the last step, we used the boundedness of $M_{\sigma}$ on $L^{p}(\sigma)$, and the bound is independent of $\sigma$. For the other direction, notice that
$$\langle f^{\nu}\rangle_{Q}^{\sigma}\leq\inf_{x\in Q}M_{\sigma}(f^{\nu})(x)=(\inf_{x\in Q}M_{\sigma,\nu}(f)(x))^{\nu}\leq(\langle M_{\sigma,\nu}(f)\rangle_{Q}^{\sigma})^{\nu},$$
where $M_{\sigma,\nu}(f):=(M_{\sigma}(f^{\nu}))^{1/\nu}$, with this observation, we have
\begin{align*}
&\sup_{\|f\|_{L^{p}(\sigma)}=1}\big\|\sum_{Q\in\mathcal{S}}\langle\sigma\rangle_{\alpha,Q}^{\nu}\langle f^{\nu}\rangle_{Q}^{\sigma}{\bf{1}}_{Q}\big\|_{L^{\frac{q}{\nu},\infty}(w)}\\
&\leq\sup_{\|f\|_{L^{p}(\sigma)}=1}\big\|\sum_{Q\in\mathcal{S}}\langle\sigma\rangle_{\alpha,Q}^{\nu}(\langle M_{\sigma,\nu}(f)\rangle_{Q}^{\sigma})^{\nu}{\bf{1}}_{Q}\big\|_{L^{\frac{q}{\nu},\infty}(w)}\\
&\leq\sup_{\|f\|_{L^{p}(\sigma)}=1}\big\|\mathcal{A}_{\alpha,\mathcal{S}}^{\nu}(\cdot\sigma)\big\|_{L^{p}(\sigma)\rightarrow L^{q,\infty}(w)}\|M_{\sigma,\nu}(f)\|_{L^{p}(\sigma)}^{\nu}\\
&\lesssim\big\|\mathcal{A}_{\alpha,\mathcal{S}}^{\nu}(\cdot\sigma)\big\|_{L^{p}(\sigma)\rightarrow L^{q,\infty}(w)},
\end{align*}
where in the last step, we use the boundedness of $M_{\sigma,\nu}$ on $L^{p}(\sigma)$ since $p>\nu$, and the bound is independent of $\sigma$. This completes the proof of Lemma $\ref{Lemma-1.1}$. $\hfill$ $\Box$

Now suppose that $B$ is the sharp constant such that
$$\big\|\sum_{Q\in\mathcal{S}}\langle\sigma\rangle_{\alpha,Q}^{\nu}\langle f^{\nu}\rangle_{Q}^{\sigma}{\bf{1}}_{Q}\big\|_{L^{\frac{q}{\nu},\infty}(w)}\leq B\|f\|_{L^{p}(\sigma)}^{\nu},$$
that is,
\begin{equation}\label{equivalence transformation}
\big\|\sum_{Q\in\mathcal{S}}\langle\sigma\rangle_{\alpha,Q}^{\nu}\langle f\rangle_{Q}^{\sigma}{\bf{1}}_{Q}\big\|_{X^{\frac{q}{\nu}}(w)}\leq B\|f\|_{L^{\frac{p}{\nu}}(\sigma)},
\end{equation}
Then
$$\|\mathcal{A}_{\alpha,\mathcal{S}}(\cdot\sigma)\|_{L^{p}(\sigma)\rightarrow X^{q}(w)}\simeq B^{\frac{1}{\nu}}.$$
Hence, we have reduced the problem to study $\eqref{equivalence transformation}$. We need the following result given by Lacey, Sawyer and Uriarte-Tuero \cite{LSU}.

\begin{prop}\label{weighted estimates for generally linear operator}
	Let $\tau=\{\tau:\,Q\in\mathcal{Q}\}$ be nonnegative constants, $w,\sigma$ be weights and define linear operators by
	$$T_{\tau}:=\sum_{Q\in\mathcal{Q}}\tau_{Q}\langle f\rangle_{Q}{\bf{1}}_{Q}.$$
	Then for $1<p\leq q<\infty$, there holds
	$$
	\|T_{\tau}(\cdot\sigma)\|_{L^{p}(\sigma)\rightarrow L^{q,\infty}(w)}\simeq\sup_{R\in\mathcal{Q}}w(R)^{-\frac{1}{q^{\prime}}}\big\|\mathop{\sum_{Q\in\mathcal{Q}}}_{Q\subset R}\tau_{Q}\langle w\rangle_{Q}{\bf{1}}_{Q}\big\|_{L^{p^{\prime}}(\sigma)}
	$$
\end{prop}
Observing that for $\eqref{equivalence transformation}$, we have
$$
\big\|\sum_{Q\in\mathcal{S}}\langle\sigma\rangle_{\alpha,Q}^{\nu}\langle f\rangle_{Q}^{\sigma}{\bf{1}}_{Q}\big\|_{L^{\frac{q}{\nu},\infty}(w)}=\|T_{\tau}(f\sigma)\|_{L^{\frac{q}{\nu},\infty}(w)}
$$
with $\tau_{Q}=\langle\sigma\rangle_{\alpha,Q}^{\nu-1}|Q|^{\frac{\alpha}{d}}$. Theorem $\ref{fractional sparse operator controlled by testing constants}$ follows immediately from Proposition $\ref{weighted estimates for generally linear operator}$.

The following proposition is weighted weak estimate for fractional maximal operator, which can found in the paper\cite{GM}.
 \begin{prop}\label{weighted estimate for fractional Maximal operator}
 	Given $1<p\leq q<\infty$, $0\leq\alpha<d$ and a pair of wights $(w,\sigma)$. Then for all measurable functions $f$,
 	$$\|M_{\alpha}(f\sigma)\|_{L^{q,\infty}(w)}\lesssim[w,\sigma]_{A_{p,q}^{\alpha}}\|f\|_{L^{p}(\sigma)}.$$
 \end{prop}

\section{Proof of the weak type bound for $1<p\leq\nu$}
We are left to prove Theorem $\ref{weak estimate for sparse operator}$ in the case that $1<p\leq\nu$. Actually, the method stem from Hyt\"{o}nen and Li \cite{HL}, they have investigated the two-weight case. Following their method, it is easy to give the off-diagonal two-weight estimate as well. For completeness, we give the deails.

{\noindent}4.1. {\bf{The case for}} $1<p\leq q<\nu$. We want to bound the following inequality,
$$
\sup_{\lambda>0}\lambda w(\{x\in\mathbb{R}^{n}:\,\mathcal{A}_{\alpha,\nu}^{\mathcal{S}}(f\sigma)>\lambda\})^{\frac{1}{q}}\lesssim [w,\sigma]_{A_{p.q}^{\alpha}}^{\frac{1}{q}}\|f\|_{L^{p}(\sigma)}.
$$
By scaling it suffices to give an uniform estimate for
$$\lambda_{0} w(\{x\in\mathbb{R}^{n}:\,\mathcal{A}_{\alpha,\nu}^{\mathcal{S}}(f\sigma)>\lambda_{0}\})^{\frac{1}{q}},$$
where $\lambda_{0}$ is some constant to be determined later. It is also free to further sparsify $\mathcal{S}$ such that
$$\big|\mathop{\bigcup_{Q^{\prime\subsetneq Q}}}_{Q^{\prime},Q\in\mathcal{S}}Q^{\prime}\big|\leq\frac{1}{4^{1-\frac{\alpha}{d}}}|Q|.$$
Now set
\begin{equation}\label{define fractional average set-1 }
\mathcal{S}_{m}:=\{Q\in\mathcal{S}:\,2^{-m-1}<\langle f\sigma\rangle_{\alpha,Q}\leq 2^{-m}\},\qquad m\geq0,
\end{equation}
and
\begin{equation}\label{define fractional average set-2 }
\mathcal{S}^{\prime}:=\{Q\in\mathcal{S}:\,\langle f\sigma\rangle_{\alpha,Q}>1\}.
\end{equation}
Then for $Q\in\mathcal{S}_{m}$, $m\geq0$, denote by $\text{ch}_{\mathcal{S}_{m}}(Q)$ the maximal subcubes of $Q$ in $\mathcal{S}_{m}$ and define
\begin{equation}\label{define E_{Q}}
E_{Q}:=Q\backslash\displaystyle\bigcup_{Q^{\prime}\in\text{ch}_{\mathcal{S}_{m}}(Q)}Q^{\prime}.
\end{equation}
Then
\begin{align}\label{related to  fractional average-1}
\langle f\sigma{\bf{1}}_{E_{Q}}\rangle_{\alpha,Q}&=\frac{1}{|Q|^{1-\frac{\alpha}{d}}}f\sigma dx-\frac{1}{|Q|^{1-\frac{\alpha}{d}}}\sum_{Q^{\prime}\in\text{ch}_{\mathcal{S}_{m}}(Q)}\int_{Q^{\prime}}f\sigma dx\nonumber\\
&=\frac{1}{|Q|^{1-\frac{\alpha}{d}}}f\sigma dx-\sum_{Q^{\prime}\in\text{ch}_{\mathcal{S}_{m}}(Q)}\left(\frac{Q^{\prime}}{|Q|}\right)^{1-\frac{\alpha}{d}}\frac{1}{|Q^{\prime}|}\int_{Q^{\prime}}f\sigma dx\\
&\geq\frac{1}{|Q|^{1-\frac{\alpha}{d}}}f\sigma dx-\frac{1}{4}2^{-l}\geq\frac{1}{2}\langle f\sigma\rangle_{\alpha,Q}\nonumber.
\end{align}
Also, we set $\mathcal{A}_{\alpha,\nu}^{\mathcal{S}_{m}}$ and $\mathcal{A}_{\alpha,\nu}^{\mathcal{S}^{\prime}}$ to be the sparse operators associated with $\mathcal{S}_{m}$ and $\mathcal{S}^{\prime}$, respectively
\begin{equation}\label{define two parts sparse operators}
(\mathcal{A}_{\alpha,\nu}^{\mathcal{S}_{m}}(f))^{\nu}:=\sum_{Q\in\mathcal{S}_{m}}\langle f\rangle_{\alpha,Q}^{\nu}{\bf{1}}_{Q}\qquad \text{and}\qquad (\mathcal{A}_{\alpha,\nu}^{\mathcal{S}^{\prime}}(f))^{\nu}:=\sum_{Q\in\mathcal{S}^{\prime}}\langle f\rangle_{\alpha,Q}^{\nu}{\bf{1}}_{Q}.
\end{equation}
Thus, it is easy to know that
\begin{equation}\label{spase divide into two parts}
\mathcal{A}_{\alpha,\nu}^{\mathcal{S}}:=\sum_{Q\in\mathcal{S}}\langle f\rangle_{\alpha,Q}^{\nu}{\bf{1}}_{Q}=\sum_{m\in\mathbb{N}}(\mathcal{A}_{\alpha,\nu}^{\mathcal{S}_{m}}(f))^{\nu}+(\mathcal{A}_{\alpha,\nu}^{\mathcal{S}^{\prime}}(f))^{\nu}.
\end{equation}
By $\eqref{define two parts sparse operators}$ and $\eqref{spase divide into two parts}$, we conclude that
\begin{align*}
&w(\{x\in\mathbb{R}^{n}:\,\mathcal{A}_{\alpha,\nu}^{\mathcal{S}}(f\sigma)>\lambda_{0}\})\\
&\leq w(\{x\in\mathbb{R}^{n}:\,\sum_{m\geq0}(\mathcal{A}_{\alpha,\nu}^{\mathcal{S}_{m}}(f))^{\nu}>\frac{\lambda_{0}^{\nu}}{2}\})+ w(\{x\in\mathbb{R}^{n}:\,(\mathcal{A}_{\alpha,\nu}^{\mathcal{S}^{\prime}}(f))^{\nu}>\frac{\lambda_{0}^{\nu}}{2}\})\\
&= w(\{x\in\mathbb{R}^{n}:\,\sum_{m\geq0}\sum_{Q\in\mathcal{S}_{m}}\langle f\sigma\rangle_{\alpha,Q}^{\nu}{\bf{1}}_{Q}>\frac{\lambda_{0}^{\nu}}{2}\})+ w(\{x\in\mathbb{R}^{n}:\,\sum_{Q\in\mathcal{S}^{\prime}}\langle f\sigma\rangle_{\alpha,Q}^{\nu}{\bf{1}}_{Q}>\frac{\lambda_{0}^{\nu}}{2}\}):=II_{1}+II_{2}.\\
\end{align*}
The second term estimation is trival. In fact, it follows immediately from Proposition $\ref{weighted estimate for fractional Maximal operator}$,
$$II_{2}\leq w\big(\bigcup_{Q\in\mathcal{S}^{\prime}}Q\big)\leq w(\{x\in\mathbb{R}^{n}:\,M_{\alpha}(f\sigma)>1\})\lesssim[w,\sigma]_{A_{p,q}^{\alpha}}\|f\|_{L^{p}(\sigma)}^{q}.$$
 Now let $\frac{\lambda_{0}^{\nu}}{2}=\displaystyle\sum_{m\geq0}2^{-\varepsilon m}$, where $\varepsilon:=(\nu-q)/2$. By $\eqref{related to  fractional average-1}$, we can estimate
\begin{align*}
II_{1}&\leq\sum_{m\geq0}w(\{x\in\mathbb{R}^{n}:\,\sum_{Q\in\mathcal{S}_{m}}\langle f\sigma\rangle_{\alpha,Q}^{\nu}{\bf{1}}_{Q}>2^{-\varepsilon m}\})\\
&\leq\sum_{m\geq0}w(\{x\in\mathbb{R}^{n}:\,\sum_{Q\in\mathcal{S}_{m}}\langle f\sigma{\bf{1}}_{Q}\rangle_{\alpha,Q}^{q}{\bf{1}}_{Q}>2^{(\nu-q)m}2^{-\varepsilon m}\})\\
&\leq\sum_{m\geq0}w(\{x\in\mathbb{R}^{n}:\,\sum_{Q\in\mathcal{S}_{m}}\langle f\sigma{\bf{1}}_{E_{Q}}\rangle_{\alpha,Q}^{q}{\bf{1}}_{Q}>2^{-q}2^{(\nu-q)m}2^{-\varepsilon m}\})\\
&\leq\sum_{m\geq0}2^{(q-\nu+\varepsilon)m+q}\int_{\mathbb{R}^{n}}\sum_{Q\in\mathcal{S}_{m}}\langle f\sigma{\bf{1}}_{E_{Q}}\rangle_{\alpha,Q}^{q}{\bf{1}}_{Q}dw\lesssim[w,\sigma]_{A_{p,q}^{\alpha}}\|f\|_{L^{p}(\sigma)}^{q}
\end{align*}
where in the last inequality we have use the following the fact
\begin{align*}
\int_{\mathbb{R}^{n}}\sum_{Q\in\mathcal{S}_{m}}\langle f\sigma{\bf{1}}_{E_{Q}}\rangle_{\alpha,Q}^{q}{\bf{1}}_{Q}dw&=\sum_{Q\in\mathcal{S}_{m}}\langle f\sigma{\bf{1}}_{E_{Q}}\rangle_{\alpha,Q}^{q}w(Q)\\
&\leq\sum_{Q\in\mathcal{S}_{m}}\left(\frac{1}{\sigma(E_{Q})^{1-\frac{1}{p}+\frac{1}{q}}}\int_{E_{Q}}f\sigma\right)^{q}|Q|^{q(\frac{\alpha}{d}-1)}w(Q)\sigma({Q})^{\frac{q}{p^{\prime}}}\sigma(E_{Q})\\
&\leq[w,\sigma]_{A_{p,q}^{\alpha}}\|f\|_{L^{p}(\sigma)}.
\end{align*}
Combining the above $II_{1}$ and $II_{2}$, we get
$$\|\mathcal{A}_{\alpha,\nu}^{\mathcal{S}}(f\sigma)\|_{L^{q,\infty}(w)}\lesssim[w,\sigma]_{A_{p,q}^{\alpha}}^{\frac{1}{q}}\|f\|_{L^{p}(\sigma)}.$$
{\noindent}4.2. {\bf{The cases for}} $p\leq q=\nu$ or $p\leq\nu<q$. We can estimate for the case by \cite[Theorem 1.1]{FH}
$$\|\mathcal{A}_{\alpha,\nu}^{\mathcal{S}}(f\sigma)\|_{L^{q,\infty}(w)}\leq\|\mathcal{A}_{\alpha,\nu}^{\mathcal{S}}(f\sigma)\|_{L^{q}(w)}\lesssim[w,\sigma]_{A_{p,q}^{\alpha}}^{\frac{1}{q}}[\sigma]_{A_{\infty}}^{\frac{1}{q}}\|f\|_{L^{p}(\sigma)}.$$
$\hfill$ $\Box$

\section{Sharpness of the weak type bounds for fractional square function}

In this section, we will show that the case for  $\nu=2$, which called  fractional square function, i.e.
\begin{equation}\label{define  fractional square function}
\mathcal{A}_{\alpha,2}^{\mathcal{S}}(f)=\big(\sum_{Q\in\mathcal{S}}\langle f\rangle_{\alpha,Q}^{2}{\bf{1}}_{Q}\big)^{\frac{1}{2}},
\end{equation}
and $p,q,\alpha$ satisfy condition $\frac{1}{q}+\frac{\alpha}{d}=\frac{1}{p}$. We only consider one weight theory estimate for $L^{p}(w^{p})\rightarrow L^{q,\infty}(w^{q})$ in here. The governing weight class is a generalization of Muckenhoupt $A_{p}$ weights, and was introduced by Muckenhoupt and Wheeden \cite{MW}.
$$[w]_{A_{p,q}}:=\sup_{Q}\left(\frac{1}{|Q|}\int_{Q}w^{q}\right)\left(\frac{1}{|Q|}\int_{Q}w^{-p^{\prime}}\right)^{\frac{q}{p^{\prime}}}<\infty.$$
Its relation to two weight characteristic is $[w^{q},w^{-p^{\prime}}]_{A_{p,q}^{\alpha}}=[w]_{A_{p,q}}$ with $\frac{1}{q}+\frac{\alpha}{d}=\frac{1}{p}$. Moreover, it is straightforward to show that the following are equivalent:
\begin{align}\label{equivalent weight}
{(\rm a)}\quad w\in A_{p,q};\qquad\qquad\qquad\qquad\text{\rm(b)}\quad w^{q}\in A_{1+\frac{q}{p^{\prime}}}\qquad\text{and} \qquad w^{-p^{\prime}}\in A_{1+\frac{p^{\prime}}{q}}.
\end{align}
We will show that the norm bound $$\|\mathcal{A}_{\alpha,2}^{\mathcal{S}}\|_{L^{p}(w^{p})\rightarrow L^{q,\infty}(w^{q})}\leq[w]_{A_{p,q}}^{\max(\frac{1}{q},\frac{1}{2}-\frac{\alpha}{d})}$$
is unimprovable. Actually, a lower bound with the exponent $\frac{1}{q}$ holds uniformly over all weights, which is the content of the next Theorem. The optimality of the exponent $\frac{1}{2}-\frac{\alpha}{d}$ is slightly more tricky, and is based on a example of a specific weight $A_{p,q}$. Also, by Theorem $\ref{weak estimate for sparse operator}$ we give the following mixed $A_{p,q}-A_{\infty}$ estimate.
\begin{coro}\label{weak estimate for square function-1}
	Let $0<\alpha<d$ and $1<p\leq q<\infty$ with $\frac{1}{q}+\frac{\alpha}{d}=\frac{1}{p}$. Then
	$$\|\mathcal{A}_{\alpha,2}^{\mathcal{S}}\|_{L^{p}(w^{p})\rightarrow L^{q,\infty}(w^{q})}\lesssim[w,\sigma]_{A_{p,q}^{\alpha}}^{\frac{1}{q}}
	\left\{
	\begin{aligned}
	&[w^{-p^{\prime}}]_{A_{\infty}}^{\frac{1}{q}}, & \quad & 2\leq q\leq\frac{2}{1-\frac{2\alpha}{d}},\\
	&[w^{q}]_{A_{\infty}}^{(\frac{1}{2}-\frac{1}{p})_{+}},  & \quad & \text{other case}.
	\end{aligned}
	\right.$$
\end{coro}
Notice that $\eqref{equivalent weight}$, we easilly know that
\begin{equation}\label{A_{p,q} tranfors to A_{p}}
[w^{q}]_{A_{1+\frac{q}{p^{\prime}}}}=[w]_{A_{p,q}}\qquad \text{and}\qquad [w^{-p^{\prime}}]_{A_{1+\frac{p^{\prime}}{q}}}=[w]_{A_{p,q}}^{\frac{p^{\prime}}{q}}.
\end{equation}
And Lerner \cite{L1} show that $[w]_{A_{\infty}}\lesssim[w]_{A_{p}}$. Hence, using this relation to Corollary $\ref{weak estimate for square function-1}$ we obtain the following pure $A_{p,q}$ estimate.
\begin{coro}\label{weak estimate for square function-2}
	Let $0<\alpha<d$ and $1<p\leq q<\infty$ with $\frac{1}{q}+\frac{\alpha}{d}=\frac{1}{p}$. Then
	$$\|\mathcal{A}_{\alpha,2}^{\mathcal{S}}\|_{L^{p}(w^{p})\rightarrow L^{q,\infty}(w^{q})}\lesssim
		\left\{
		\begin{aligned}
		&[w]_{A_{p,q}}^{\frac{p^{\prime}}{q}(1-\frac{\alpha}{d})}, & \quad & 2\leq q\leq\frac{2}{1-\frac{2\alpha}{d}},\\
		&[w]_{A_{p,q}}^{(\frac{1}{2}-\frac{\alpha}{d})},  & \quad & \frac{2}{1-\frac{2\alpha}{d}}<q<\infty,\\
		&[w]_{A_{p,q}}^{\frac{1}{q}},  & \quad & 1\leq q<2.
		\end{aligned}
		\right.
	$$
\end{coro}
However, the exponent $\frac{p^{\prime}}{q}(1-\frac{\alpha}{d})$ is not optimal of  the case for $2\leq q\leq\frac{2}{1-\frac{2\alpha}{d}}$. and the best  exponent $\frac{1}{q}$ will appear following estimate.
For generally, we consider case for $\nu\geq1$, and we are concerned with the weak-type bounds, which have cases of $\nu\leq q\leq\frac{\nu}{1-\frac{\nu\alpha}{d}}$, which a $(\log_{1}[w^{q}]_{A_{\infty}})^{\frac{1}{\nu}}$ appears in the sharp estimate.
\begin{theorem}\label{one-weight sharp estimates for sparse operators}
	Let $\nu\geq1$, $0\leq\alpha<d$ and $1\leq p\leq q<\infty$ with $\frac{1}{p}+\frac{\alpha}{d}=\frac{1}{q}$, there holds for any weight $w\in A_{p,q}$
	$$\|\mathcal{A}_{\alpha,\nu}^{\mathcal{S}}(f)\|_{L^{q,\infty}(w^{q})}\lesssim [w]_{A_{p,q}}^{\max(\frac{1}{q},\frac{1}{\nu}-\frac{\alpha}{d})}\phi([w^{q}]_{A_{\infty}})\|wf\|_{L^{p}},$$
	where
	$$
	\phi([w^{q}]_{A_{\infty}})=
	\left\{
	\begin{aligned}
	&(\log_{1}[w^{q}]_{A_{\infty}})^{\frac{1}{\nu}}, & \quad & \nu\leq q\leq\frac{\nu}{1-\frac{\nu\alpha}{d}};\\
	&1, & \quad & \text{other case}. \\
	\end{aligned}
	\right.
	$$
	and $\log_{1}(x)=1+\log_{+}(x)$.
\end{theorem}

As a Corollary of Theorem $\ref{one-weight sharp estimates for sparse operators}$, the following result of fractional square function is sharp.

\begin{coro}\label{one-weight sharp estimates for fractional square functions}
	Let $0\leq\alpha<d$ and $1\leq p\leq q<\infty$ with $\frac{1}{p}+\frac{\alpha}{d}=\frac{1}{q}$, there holds for any weight $w\in A_{p,q}$
	$$\|\mathcal{A}_{\alpha,2}^{\mathcal{S}}(f)\|_{L^{q,\infty}(w^{q})}\lesssim [w]_{A_{p,q}}^{\max(\frac{1}{q},\frac{1}{2}-\frac{\alpha}{d})}\phi_{1}([w^{q}]_{A_{\infty}})\|wf\|_{L^{p}},$$
	where
	$$
	\phi_{1}([w^{q}]_{A_{\infty}})=
	\left\{
	\begin{aligned}
	&(\log_{1}[w^{q}]_{A_{\infty}})^{\frac{1}{2}}, & \quad & 2\leq q\leq\frac{2}{1-\frac{2\alpha}{d}};\\
	&1, & \quad & \text{other case}. \\
	\end{aligned}
	\right.
	$$
\end{coro}

A basic tool for us is the following classical reverse H\"{o}lder's inequality with
optimal bound, which can be found in \cite{HP}.

\begin{prop}\label{proposition for reverse H?lder estimate}
	There is a dimensional constant $c>0$ such that for $w\in A_{\infty}$, and $r(w)=1+c[w]_{A_{\infty}}$, there holds
	\begin{equation}\label{reverse H?lder estimate}
	\langle w^{r(w)}\rangle_{Q}^{\frac{1}{r(w)}}\leq2\langle w\rangle_{Q}, \qquad\qquad  Q \,\,\,\text{\rm a cube.}
	\end{equation}
\end{prop}

We also need the following off-diagonal extrapolation given by Duoandikoetxra \cite{D1}.
\begin{prop}\label{ off-diagonal extrapolation}
	Let $1\leq p_{0}<\infty$ and $0<q_{0}<\infty$. Assume that that for some family of nonnegative couples $(f,g)$ and for all $w\in A_{p_{0},q_{0}}$ we have
	\begin{align*}
	\|wg\|_{L^{q_{0}}}\leq CN([w]_{A_{p_{0},q_{0}}})\|wf\|_{L^{p_{0}}},
	\end{align*}
	where $N$ is an increasing function and the constant $C$ does not depend on $w$. Set $\gamma=\frac{1}{q_{0}}+\frac{1}{p_{0}^{\prime}}$. Then for $1<p<\infty$ and $0<q<\infty$, such that
	$$\frac{1}{q}-\frac{1}{p}=\frac{1}{q_{0}}-\frac{1}{q_{0}},$$
	and all $w\in A_{p,q}$ we have
	$$\|wg\|_{L^{p}}\leq CK(w)\|wf\|_{L^{p}},$$
	where
		$$
		K(w)=
		\left\{
		\begin{aligned}
		&N([w]_{A_{p,q}}(2\|M\|_{L^{\gamma q}(w^{q})})^{\gamma(q-q_{0})}), & \quad & q<q_{0}; \\
		&N([w]_{A_{p,q}}^{\frac{\gamma q_{0}-1}{\gamma q-1}}(2\|M\|_{L^{\gamma p^{\prime}}(w^{-p^{\prime}})})^{\frac{\gamma(q-q_{0})}{\gamma q-1}}), & \quad & q>q_{0}.\\
		\end{aligned}
		\right.
		$$
		In particular, $K(w)\leq C_{1}N(C_{2}[w]_{A_{p,q}}^{\max(1,\frac{q_{0}p^{\prime}}{qp_{0}^{\prime}})})$ for $w\in A_{p,q}$.
\end{prop}

The following estimate based on Domingo-Salazar, Lacey, and Rey \cite{DLR}.
\begin{theorem}\label{sharp weighted for fractional suqare function-1}
	Let $\nu\geq1$, $0\leq\alpha<d$ and $1\leq p\leq q<\infty$ with $\frac{1}{p}+\frac{\alpha}{d}=\frac{1}{q}$, there holds for any weight $w\in A_{p,q}$
	\begin{equation}\label{inequality of sharp weak weighted for sparse}
	\|\mathcal{A}_{\alpha,\nu}^{\mathcal{S}}(f)\|_{L^{q,\infty}(w^{q})}\lesssim[w]_{A_{p,q}}^{\max(\frac{1}{q},\frac{1}{\nu}-\frac{\alpha}{d})}\phi_{2}([w^{q}]_{A_{\infty}})\|wf\|_{L^{p}},
	\end{equation}
	where
	$$
	\phi_{2}([w^{q}]_{A_{\infty}})=
	\left\{
	\begin{aligned}
	&1, & \quad & 1\leq q<\nu; \\
	&(\log_{1}[w^{q}]_{A_{\infty}})^{\frac{1}{\nu}}, & \quad &\nu\leq q<\infty.\\
	\end{aligned}
	\right.
	$$
\end{theorem}
 {\noindent}Theorem $\ref{one-weight sharp estimates for sparse operators}$ follows immediately from Theorems $\ref{weak estimate for sparse operator}$ and $\ref{sharp weighted for fractional suqare function-1}$.

In order to prove Theorem $\ref{sharp weighted for fractional suqare function-1}$, we need following estimate.
\begin{lemma}\label{lemma of fractonal square function-1}
	Let $\nu\geq1$, $q\geq \nu$, $0\leq\alpha<d$ and $1\leq p\leq q<\infty$ with $\frac{1}{q}+\frac{\alpha}{d}=\frac{1}{p}$, then
	$$
	\|\mathcal{A}_{\alpha,\nu}^{\mathcal{S}_{m}}\|_{L^{q,\infty}(w^{q})}\lesssim [w]_{A_{p,q}}^{\frac{1}{\nu}-\frac{\alpha}{d}}\|wf\|_{L^{p}}.
	$$
	where  and given $0<m<\log_{1}[w^{q}]_{A_{\infty}}$.
\end{lemma}

{\noindent}$Proof$. We only need to prove the case for $q=\frac{\nu}{1-\frac{\nu\alpha}{d}}$, by off-diagonal extrapolation in Proposition $\ref{ off-diagonal extrapolation}$, it yields the cases for $\nu\leq q<\frac{\nu}{1-\frac{\nu\alpha}{d}}$ and $\frac{\nu}{1-\frac{\nu\alpha}{d}}<q<\infty$. By Minkowski's inequality and  $\eqref{related to  fractional average-1}$, we can estimate
\begin{align*}
\left(\int_{\mathbb{R}^{n}}\big(\sum_{Q\in\mathcal{S}_{m}}\langle f\rangle_{\alpha,Q}^{\nu}{\bf{1}}_{Q}\big)^{\frac{q}{\nu}}\upsilon\right)^{\frac{1}{q}}&\leq\left(\sum_{Q\in\mathcal{S}_{m}}\left(\int_{\mathbb{R}^{n}}\langle f\rangle_{\alpha,Q}^{q}{\bf{1}}_{Q}\upsilon\right)^{\frac{\nu}{q}}\right)^{\frac{1}{\nu}}=\left(\sum_{Q\in\mathcal{S}_{m}}\langle f\rangle_{\alpha,Q}^{\nu}\upsilon^{\frac{\nu}{q}}(Q)\right)^{\frac{1}{\nu}}\\
&\lesssim\left(\sum_{Q\in\mathcal{S}_{m}}\langle f{\bf{1}}_{E_{m}(Q)}\rangle_{\alpha,Q}^{\nu}\upsilon^{\frac{\nu}{q}}(Q)\right)^{\frac{1}{\nu}}\\
&\leq\left(\sum_{Q\in\mathcal{S}_{m}}\langle f^{p}{\bf{1}}_{E_{m}(Q)}w^{p}\rangle_{\alpha,Q}\langle \sigma^{p^{\prime}}\rangle_{\alpha,Q}^{\frac{p}{p^{\prime}}}\upsilon^{\frac{p}{q}}(Q)\right)^{\frac{1}{p}}\\
&\leq[w]_{A_{q,\nu}}^{\frac{1}{q}}\left(\sum_{Q\in\mathcal{S}_{m}}\int_{E_{m}(Q)}f^{p}{\bf{1}}_{E_{m}(Q)}w^{p}\right)^{\frac{1}{p}}\leq[w]_{A_{q,\nu}}^{\frac{1}{\nu}-\frac{\alpha}{d}}\|wf\|_{L^{p}},
\end{align*}
where $p=\nu$ in the above estimate. $\hfill$ $\Box$

The good property of Lebesgue measure appear in the paper \cite{DLR}.
\begin{prop}\label{good property of Lebesgue measure}
	Let any $\lambda>0$, $\mathcal{S}_{m}$ defined as $\eqref{define fractional average set-1 }$ and $b=\sum_{Q^{\prime}\in\mathcal{S}_{m}}{\bf{1}}_{Q^{\prime}}$, then we have that for any dyadic cube $Q\in\mathcal{S}_{m}$ 
	$$|\{x\in Q:\,b(x)>\lambda\}|\lesssim \exp(-c\lambda)|Q|.$$
\end{prop}

For $\log_{1}[w^{q}]_{A_{\infty}}\leq m<\infty$, we also have following estimate.

\begin{lemma}\label{lemma of fractonal square function-2}
	Let $\upsilon$ denote the weight $w^{q}$, for all integers $m_{0}>0$, then
	\begin{equation}\label{weak estimate for first part sparse operator}
	\upsilon\big(\sum_{m=m_{0}}^{\infty}(\mathcal{A}_{\alpha,\nu}^{\mathcal{S}_{m}}(f))^{\nu}>1\big)\lesssim[w]_{A_{p,q}}\left(\frac{[w]_{A_{\infty}}}{2^{m_{0}}}\right)^{q}\|wf\|_{L^{p}}^{q}.
	\end{equation}
\end{lemma}
{\noindent}$Proof$. Define 
$$
\mathcal{S}_{m}^{*}:=\{Q\,\,\text{maximal\,\,s.t.}\,\,Q\in\mathcal{S}_{m}\}\quad\text{and}\quad B_{m}:=\bigcup\{Q:\,Q\in\mathcal{S}_{m}^{*}\}.
$$
By the definitions of $\mathcal{S}_{m}$ and $(\mathcal{A}_{\alpha,\nu}^{\mathcal{S}_{m}}(f))^{\nu}$, we can write $(\mathcal{A}_{\alpha,\nu}^{\mathcal{S}_{m}}(f))^{\nu}$ as $2^{-\nu
	m}b_{m}$, where
$$b_{m}\leq\sum_{Q\in\mathcal{S}_{m}}{\bf{1}}_{Q}\quad\text{and}\quad \text{supp}(b_{m})\subset B_{m}.$$
For any dyadic cube $Q\in\mathcal{S}_{m}$, by Proposition $\ref{good property of Lebesgue measure}$, we know that the function $b_{m}$ is locally exponentially integrable. By the sharp weak-type estimate for the  fractional maximal function \cite{LMPT}, we know that
$$\upsilon(B_{m})\lesssim2^{qm}[w]_{A_{p,q}}\|wf\|_{L^{p}}^{q}.$$
The left hand side of $\eqref{weak estimate for first part sparse operator}$ can be estimated as
\begin{align*}
\upsilon\big(\sum_{m=m_{0}}^{\infty}(\mathcal{A}_{\alpha,\nu}^{\mathcal{S}_{m}}(f))^{\nu}>1\big)&=\upsilon\big(\sum_{m=m_{0}}^{\infty}2^{-\nu m}b_{m}>\sum_{m=m_{0}}^{\infty}2^{m_{0}-m-1}\big)\\
&\leq\sum_{m=m_{0}}^{\infty}\upsilon(b_{m}>2^{m_{0}+(\nu-1)m-1}).
\end{align*}
Taking 
$$\beta(Q):=\{x\in Q:\,b_{m}(x)>2^{m_{0}+(\nu-1)m-1}\}$$ 
for any dyadic cube $Q\in\mathcal{S}_{m}^{*}$, by the definition of $\mathcal{S}_{m}^{*}$ and Proposition $\ref{good property of Lebesgue measure}$, we show that
$$|\beta(Q)|\lesssim\exp(-c2^{m_{0}+(\nu-1)m})|Q|.$$
Using the $A_{\infty}$ property for $A_{1+\frac{q}{p^{\prime}}}$ weights with $\upsilon$-measure and Proposition $\ref{reverse H?lder estimate}$, there holds
\begin{align*}
\upsilon(\beta(Q))&=\langle \upsilon{\bf{1}}_{\beta(Q)}\rangle_{Q}|Q|\leq\langle {\bf{1}}_{\beta(Q)}\rangle_{Q}^{(\frac{1}{r(\upsilon)})^{\prime}}\langle \upsilon^{r(\upsilon)}\rangle_{Q}^{\frac{1}{r(\upsilon)}}|Q|\\
&\lesssim \left[\frac{|\beta(Q)|}{|Q|}\right]^{(c[\upsilon]_{A_{\infty}})^{-1}}\upsilon(Q)\lesssim\upsilon(Q)\exp\left(-c\frac{2^{m_{0}+(\nu-1)m}}{[\upsilon]_{A_{\infty}}}\right),
\end{align*}
where $r(\upsilon)$ as in $\eqref{reverse H?lder estimate}$.

Summing over the disjoint cubes in $\mathcal{S}_{m}^{*}$, we obtain
\begin{equation}\label{Summing over the disjoint cubes}
\upsilon\big(\sum_{m=m_{0}}^{\infty}(\mathcal{A}_{\alpha,\nu}^{\mathcal{S}_{m}}(f))^{\nu}>1\big)\lesssim[w]_{A_{p,q}}\|wf\|_{L^{p}}^{q}\sum_{m=m_{0}}^{\infty}2^{mq}\exp\left(-c\frac{2^{m_{0}+(\nu-1)m}}{[\upsilon]_{A_{\infty}}}\right).
\end{equation}
The sum in the right hand side of $\eqref{Summing over the disjoint cubes}$, we can be controlled by
\begin{align}\label{the last sum}
\sum_{m=m_{0}}^{\infty}2^{mq}\exp\left(-c\frac{2^{m_{0}+(\nu-1)m}}{[\upsilon]_{A_{\infty}}}\right)&\leq\int_{m_{0}}^{\infty}2^{qx}\exp\left(-c\frac{2^{m_{0}+(\nu-1)x}}{[\upsilon]_{A_{\infty}}}\right)dx\nonumber\\
&\approx\int_{2^{(\nu-1)m_{0}}}^{\infty}y^{q}\exp\left(-c\frac{2^{m_{0}}}{[\upsilon]_{\infty}}y\right)\frac{dy}{y}\nonumber\\
&=\left(\frac{[\upsilon]_{\infty}}{2^{m_{0}}}\right)^{q}\int_{\frac{2^{\nu m_{0}}}{[\upsilon]_{\infty}}}^{\infty}y^{q}e^{-cy}\frac{dy}{y}\lesssim
\left(\frac{[\upsilon]_{\infty}}{2^{m_{0}}}\right)^{q}.
\end{align}
Combining $\eqref{Summing over the disjoint cubes}$ and $\eqref{the last sum}$, we obtain the desired result. This completes the proof Lemma $\ref{lemma of fractonal square function-2}$. $\hfill$ $\Box$

$Proof$ $of$ $Theorem$ $\ref{sharp weighted for fractional suqare function-1}$. The case for $1\leq q<\nu$ is easy and contained in Theorem $\ref{weak estimate for sparse operator}$,as so only our attention on the case for $q\geq\nu$. By scaling the left hand side of $\eqref{inequality of sharp weak weighted for sparse}$ suffices to estimate 
\begin{equation}\label{suffices to estimate for sparse operator}
\lambda^{q}\upsilon(\{x\in\mathbb{R}^{n}:\,\mathcal{A}_{\alpha,\nu}^{\mathcal{S}}>\lambda\}).
\end{equation}
Now, we assume that $\lambda=3^{\frac{1}{\nu}}$, $\|f\|_{L^{p}(w^{p})}=1$ and notice that $\eqref{spase divide into two parts}$, the $\eqref{suffices to estimate for sparse operator}$ can be estimated as
\begin{align*}
\upsilon((\mathcal{A}_{\alpha,\nu}^{\mathcal{S}}(f))^{\nu}>3)\leq\upsilon((\mathcal{A}_{\alpha,\nu}^{\mathcal{S}^{\prime}}(f))^{\nu}>1)+\upsilon\bigg(\sum_{m=0}^{m_{0}-1}(\mathcal{A}_{\alpha,\nu}^{\mathcal{S}_{m}}(f))^{\nu}>1\bigg)+\upsilon\bigg(\sum_{m=m_{0}}^{\infty}(\mathcal{A}_{\alpha,\nu}^{\mathcal{S}_{m}}(f))^{\nu}>1\bigg).
\end{align*}
By the sharp weak-type estimate for the fractional maximal function \cite{LMPT}, the first term to arrive at the bound
\begin{equation}\label{first term estimate}
\upsilon((\mathcal{A}_{\alpha,\nu}^{\mathcal{S}^{\prime}}(f))^{\nu}>1)\lesssim [w]_{A_{p,q}}^{\frac{1}{q}}.
\end{equation}
By Chebysheff's inequality and  Minkowski's inequality for $q\geq\nu$, the second term from Lemma $\ref{lemma of fractonal square function-1}$
\begin{align}\label{second term estimate}
\upsilon\bigg(\sum_{m=0}^{m_{0}-1}(\mathcal{A}_{\alpha,\nu}^{\mathcal{S}_{m}}(f))^{\nu}>1\bigg)&\leq \bigg\|\sum_{m=0}^{m_{0}-1}(\mathcal{A}_{\alpha,\nu}^{\mathcal{S}_{m}}(f))^{\nu}\bigg\|_{L^{\frac{q}{\nu}}(\upsilon)}^{\frac{q}{\nu}}\nonumber\\
&\leq\bigg(\sum_{m=0}^{m_{0}-1}\|(\mathcal{A}_{\alpha,\nu}^{\mathcal{S}_{m}}(f))^{\nu}\|_{L^{\frac{q}{\nu}}(\upsilon)}\bigg)^{\frac{q}{\nu}}\nonumber\\
&=\bigg(\sum_{m=0}^{m_{0}-1}\|\mathcal{A}_{\alpha,\nu}^{\mathcal{S}_{m}}(f)\|_{L^{q}(w^{q})}^{\nu}\bigg)^{\frac{q}{\nu}}\lesssim(m_{0}[w]_{A_{p,q}}^{\frac{1}{\nu}-\frac{\alpha}{d}})^{\frac{q}{\nu}}.
\end{align}
By Lemma $\ref{lemma of fractonal square function-2}$, the third term  can be estimate as
\begin{equation}\label{third term estimate}
\upsilon\bigg(\sum_{m=m_{0}}^{\infty}(\mathcal{A}_{\alpha,\nu}^{\mathcal{S}_{m}}(f))^{\nu}>1\bigg)\lesssim[w]_{A_{p,q}}\left(\frac{[w^{q}]_{\infty}}{2^{m_{0}}}\right)^{q}.
\end{equation}
Combining $\eqref{first term estimate}$, $\eqref{second term estimate}$ and $\eqref{third term estimate}$, we get
$$\|\mathcal{A}_{\alpha,\nu}^{\mathcal{S}}\|_{L^{q,\infty}(w^{q})}\lesssim [w]_{A_{p,q}}^{\frac{1}{q}}+ m_{0}^{\frac{1}{\nu}}[w]_{A_{p,q}}^{\frac{1}{\nu}-\frac{\alpha}{d}}+[w]_{A_{p,q}}^{\frac{1}{q}}[w^{q}]_{A_{\infty}}2^{-m_{0}}\approx[w]_{A_{p,q}}^{\max(\frac{1}{q},\frac{1}{\nu}-\frac{\alpha}{d})}(\log_{1}[w^{q}]_{A_{\infty}})^{\frac{1}{\nu}},$$
due to $m_{0}\approx\log_{1}[w^{q}]_{A_{\infty}}$. This finishs the proof Theorem $\ref{sharp weighted for fractional suqare function-1}$.
$\hfill$ $\Box$

However, this is not the end of the story; we can prove even more. Here we present our full statement of the main theorem. This estimate is sharp in the following sense.

  \begin{theorem}\label{sharp weak estimate for square function-1}
  	For any weight $w$, we have
  	$$\|\mathcal{A}_{\alpha,\nu}^{\mathcal{S}}\|_{L^{p}(w^{p})\rightarrow L^{q,\infty}(w^{q})}\geq[w]_{A_{p,q}}^{\frac{1}{q}}.$$
  \end{theorem}
  $Proof$. Let $\upsilon$ denote the weight $w^{q}$ and  consider $f=|f|\chi_{Q}$, then we obtain for $Q\in\mathcal{S}$
  $$\mathcal{A}_{\alpha,2}^{\mathcal{S}}(f)\geq\langle |f|\rangle_{\alpha,Q}.$$
  Taking $N:=\|\mathcal{A}_{\alpha,\nu}^{\mathcal{S}}(f)\|_{L^{p}(w^{p})\rightarrow L^{q,\infty}(\upsilon)}$, by the inequality of norm $\mathcal{A}_{\alpha,\nu}^{\mathcal{S}}(f)$, we have
  $$N\|f\|_{L^{p}(w^{p})}\geq\|\mathcal{A}_{\alpha,\nu}^{\mathcal{S}}(f)\|_{L^{q,\infty}(\upsilon)}\geq\|\langle |f|\rangle_{\alpha,Q}\|_{L^{q,\infty}(\upsilon)}=\frac{\upsilon(Q)^{\frac{1}{q}}}{|Q|^{1-\frac{\alpha}{d}}}\int_{Q}|f|=\frac{\upsilon(Q)^{\frac{1}{q}}}{|Q|^{1-\frac{\alpha}{d}}}\int_{Q}|f|w^{-p}w^{p}$$
  for all positive functions $|f|$ on $Q$. By the converse to H\"{o}lder's inequality, this shows that
  $$N\geq\frac{\upsilon(Q)^{\frac{1}{q}}}{|Q|^{1-\frac{\alpha}{d}}}\|w^{-p}\|_{L^{p^{\prime}}(w^{p})}=\frac{\upsilon(Q)^{\frac{1}{q}}\sigma(Q)^{\frac{1}{p^{\prime}}}}{|Q|^{1-\frac{\alpha}{d}}},$$
  and taking the supremuum over all $Q$ proves this theorem.
  $\hfill$ $\Box$

  \begin{theorem}\label{sharp weak estimate for square function-2}
  	Let $\nu\geq1$, $0\leq\alpha<d$ and $1\leq p\leq q<\infty$ with $\frac{1}{q}+\frac{\alpha}{d}=\frac{1}{p}$. If $\Phi$ be an increasing function such that
  	$$\|\mathcal{A}_{\alpha,\nu}^{\mathcal{S}}\|_{L^{p}(w^{p})\rightarrow L^{q,\infty}(w^{q})}\leq\Phi([w]_{A_{p,q}})$$
  	for all $w\in A_{p,q}$, then $\Phi(t)\gtrsim ct^{\frac{1}{\nu}-\frac{\alpha}{d}}$.
  \end{theorem}

 Lacey and Scurry \cite{LS} show that this in class of power functions, namely, they proved that there cannot be a bound the form $\Phi(t)=t^{\frac{1}{2}-\eta}$ for $\eta>0$. We will extend their methods to general case.

  {\noindent}$Proof$. We will consider two cases to prove this theorem: $\nu>1$ and $\nu=1$.

  {\noindent}Case 1: $\nu>1$. Following the same arguments as that in \cite{HL,LS}, the assumption implies
  \begin{equation}\label{Khintchine inequality}
  \big\|\big(\sum_{Q}\langle a_{Q}\cdot w^{q}\rangle_{\alpha,Q}^{\nu}{\bf{1}}_{Q}\big)^{\frac{1}{\nu}}\big\|_{L^{p^{\prime}}(w^{-p^{\prime}})}\lesssim\Phi([w]_{A_{p,q}})\big\|\big(\sum_{Q}a_{Q}^{\nu}\big)^{\frac{1}{\nu}}\big\|_{L^{q^{\prime,1}}(w^{q})}
  \end{equation}
  for all sequences of measurable functions $a_{Q}$.  For $\vartheta\in(0,1)$, we consider $w(x)=|x|^{\frac{\vartheta-1}{q}}$ and a sequence of functions
  $$a_{[0,2^{-k})}(x):=a_{k}(x):=\vartheta^{\frac{1}{\nu-1}-\frac{1}{\nu}}\sum_{j=k+1}^{\infty}2^{-\vartheta(j-k)}{\bf{1}}_{[2^{-j},2^{-j+1})}(x),\qquad k\in\mathbb{N}.$$
  Then it is easy to check that
  $$[w]_{A_{p,q}}=[w^{q}]_{A_{1+\frac{q}{p^{\prime}}}}\simeq\vartheta^{-1}\qquad \text{and}\qquad \displaystyle\sum_{k}a_{k}^{\nu}(x)\lesssim\vartheta^{\frac{\nu}{\nu-1}-2}{\bf{1}}_{[0,1]}.$$ In fact, we choose $I_{k}=[0,2^{-k}]$ and $x\in(2^{-(l+1)},2^{-l}]$ with $l\in\mathbb{N}_{0}$ such that
  $$a_{k}(x)\simeq \vartheta^{\frac{1}{\nu-1}-\frac{1}{\nu}}|I_{k}|^{-\vartheta}|x|^{\vartheta}{\bf{1}}_{I_{k}}(x).$$
  A simple calculation shows that
  \begin{align*}
  \sum_{k=0}^{\infty}a_{k}^{\nu}(x)&=\vartheta^{\frac{\nu}{\nu-1}-1}|x|^{\nu\vartheta}\sum_{k=0}^{\infty}|I_{k}|^{-\nu\vartheta}{\bf{1}}_{I_{k}}(x)=\vartheta^{\frac{\nu}{\nu-1}-1}|x|^{\nu\vartheta}\sum_{k=0}^{l}(2^{\nu\vartheta})^{k}\\
  &=\vartheta^{\frac{\nu}{\nu-1}-1}|x|^{\nu\vartheta}\frac{2^{\nu(l+1)\vartheta}-1}{2^{\nu\vartheta}-1}\lesssim\vartheta^{\frac{\nu}{\nu-1}-2}|x|^{\nu\vartheta}2^{\nu l\vartheta}\lesssim\vartheta^{\frac{\nu}{\nu-1}-2}{\bf{1}}_{[0,1]}.
  \end{align*}
  This directly for the right hand side of $\eqref{Khintchine inequality}$
  \begin{align*}
  \big\|\big(\sum_{k=1}^{\infty}a_{k}(x)^{\nu}\big)^{\frac{1}{\nu}}\big\|_{L^{q^{\prime},1}(w^{q})}&\lesssim q^{\prime}\int_{0}^{\infty}\left(\int_{\{x\in[0,1]:\,c\vartheta^{\frac{\nu}{\nu-1}-2}>s\}}|x|^{\vartheta-1}dx\right)^{\frac{1}{q^{\prime}}}ds\\
  &\leq\int_{0}^{c\vartheta^{\frac{\nu}{\nu-1}-2}}\left(\int_{0}^{1}|x|^{\vartheta-1}dx\right)^{\frac{1}{q^{\prime}}}ds\simeq\vartheta^{\frac{\nu}{\nu-1}-2}\vartheta^{-\frac{1}{q^{\prime}}}.
  \end{align*}
  On the other hand, the left hand side of $\eqref{Khintchine inequality}$ can be estimated as
  \begin{align*}
  \langle a_{k}\cdot w^{q}\rangle_{\alpha,[0,2^{-k})}&\simeq\vartheta^{\frac{1}{\nu-1}-\frac{1}{\nu}}2^{k(1-\frac{\alpha}{d})}\sum_{j=k+1}^{\infty}2^{-\vartheta(j-k)}2^{-\vartheta j}\simeq\vartheta^{\frac{1}{\nu-1}-\frac{1}{\nu}-1}2^{k(1-\frac{\alpha}{d}-\vartheta)},
  \end{align*}
  It follows that

  \begin{align*}
  \int_{[0,1]}\big(\sum_{k=1}^{\infty}\langle a_{k}\cdot w^{q}\rangle_{\alpha,[0,2^{-k})}^{\nu}{\bf{1}}_{[0,2^{-k})}\big)^{\frac{p^{\prime}}{\nu}}w^{-p^{\prime}}&\simeq\vartheta^{\frac{p^{\prime}}{\nu-1}-\frac{p^{\prime}}{\nu}-p^{\prime}}\int_{0}^{1}|x|^{(\vartheta-(1-\frac{\alpha}{d})p^{\prime}}|x|^{-\frac{(\vartheta-1)p^{\prime}}{q}}dx\\
  &=\vartheta^{\frac{p^{\prime}}{\nu-1}-\frac{p^{\prime}}{\nu}-p^{\prime}}\int_{0}^{1}|x|^{\frac{\vartheta p^{\prime}}{q^{\prime}}-1}dx=\frac{q^{\prime}}{p^{\prime}}\vartheta^{\frac{p^{\prime}}{\nu-1}-\frac{p^{\prime}}{\nu}-p^{\prime}-1}.
  \end{align*}
  By assumption, this implies
  $$\vartheta^{\frac{1}{\nu-1}-\frac{1}{p^{\prime}}-\frac{1}{\nu}-1}\lesssim\Phi([w]_{A_{p,q}})\vartheta^{\frac{\nu}{\nu-1}-2}\vartheta^{-\frac{1}{q^{\prime}}}\lesssim\Phi(c\vartheta^{-1})\vartheta^{\frac{\nu}{\nu-1}-2}\vartheta^{-\frac{1}{q^{\prime}}}.$$
  Hence, we show that $\Phi(t)\gtrsim t^{\frac{1}{\nu}-\frac{\alpha}{d}}$, this finishes the Case 1 of the estimate.

  {\noindent}Case 2: $\nu=1$. This case upper bound follows from \cite{LMPT}, and we show that
  \begin{equation}\label{prove case 2}
  \|\mathcal{A}_{\alpha,1}^{\mathcal{S}}(f)\|_{L^{q,\infty}(w^{q})}\lesssim\Phi([w]_{A_{p,q}})\|wf\|_{L^{p}}
  \end{equation}
  holds for $\Phi(t)\geq ct^{1-\frac{\alpha}{d}}$.

  By $\eqref{A_{p,q} tranfors to A_{p}}$, we show that
  \begin{equation}\label{prove case 2_1}
  \|\mathcal{A}_{\alpha,1}^{\mathcal{S}}(f)\|_{L^{q,\infty}(w^{q})}\lesssim\Phi([w^{q}]_{A_{1+{q}/{p^{\prime}}}})\|wf\|_{L^{p}},
  \end{equation}
  and if we let $u=w^{q}$, then
  \begin{equation}\label{prove case 2_2}
  \|\mathcal{A}_{\alpha,1}^{\mathcal{S}}(f)\|_{L^{q,\infty}(u)}\lesssim\Phi([u]_{A_{1+{q}/{p^{\prime}}}})\|f\|_{L^{p}(u^{p/q})}.
  \end{equation}
  Assume now that $u\in A_{1}$, then $\eqref{prove case 2_2}$ it yields that
  \begin{equation}\label{prove case 2_3}
  \|\mathcal{A}_{\alpha,1}^{\mathcal{S}}(f)\|_{L^{q,\infty}(u)}\lesssim\Phi([u]_{A_{1}})\|f\|_{L^{p}(u^{p/q})}.
  \end{equation}
  Since $\frac{p}{q}=1-\frac{p\alpha}{d}$, this is equivalent to
  \begin{equation}\label{prove case 2_4}
  \|\mathcal{A}_{\alpha,1}^{\mathcal{S}}(u^{\frac{\alpha}{d}}f)\|_{L^{q,\infty}(u)}\lesssim\Phi([u]_{A_{1}})\|f\|_{L^{p}(u)}.
  \end{equation}
  We now prove that $\eqref{prove case 2_4}$ holds for $\Phi(t)\geq ct^{1-\frac{\alpha}{d}}$. Let
  $$u(x)=|x|^{\vartheta-n}$$
  with $0<\vartheta<1$. Then standard computations shows that
  \begin{equation}\label{A_1 constant}
  [u]_{A_{1}}\simeq{\vartheta}^{-1}.
  \end{equation}
  Consider the function $f=\chi_{B}$ where $B$ is the unit ball, we can compute its norm to be
  \begin{equation}\label{L^{p} value}
  \|f\|_{L^{p}(u)}=u(B)^{\frac{1}{p}}\simeq\vartheta^{-\frac{1}{p}}.
  \end{equation}
  By sparse domination formula, we know there exists a sparse family $\mathcal{S}$ such that
  \begin{equation}\label{sparse operator relation to fractional integral operator}
  \mathcal{A}_{\alpha,1}^{\mathcal{S}}(|f|)(x)\gtrsim |I_{\alpha}f(x)|,
  \end{equation}
 where $I_{\alpha}$ is defined by $\eqref{define fractional integral operator}$. Let $0 < x_{\vartheta} < 1$ be a parameter whose value will be chosen soon. By $\eqref{sparse operator relation to fractional integral operator}$, we have that
  \begin{align*}
  \|\mathcal{A}_{\alpha,1}^{\mathcal{S}}(u^{\frac{\alpha}{d}}f)\|_{L^{q,\infty}(u)}&\gtrsim\|I_{\alpha}u^{\frac{\alpha}{d}}f\|_{L^{q,\infty}(u)}\\
  &\geq\sup_{\lambda>0}\left(u\{|x|<x_{\vartheta}:\,\int_{B}\frac{|y|^{(\vartheta-1)\alpha/d}}{|x-y|^{1-\alpha/d}}dx>\lambda\}\right)^{\frac{1}{q}}\\
  &\geq\sup_{\lambda>0}\left(u\{|x|<x_{\vartheta}:\,\int_{B\backslash B(0,|x|)}\frac{|y|^{(\vartheta-1)\alpha/d}}{|x-y|^{1-\alpha/d}}dx>\lambda\}\right)^{\frac{1}{q}}\\
  &\geq\sup_{\lambda>0}\left(u\{|x|<x_{\vartheta}:\,\int_{B\backslash B(0,|x|)}\frac{|y|^{(\vartheta-1)\alpha/d}}{(2|y|)^{1-\alpha/d}}dx>\lambda\}\right)^{\frac{1}{q}}\\
  &=\sup_{\lambda>0}\left(u\{|x|<x_{\vartheta}:\,\frac{c_{\alpha,d}}{\vartheta}(1-|x|^{\vartheta\alpha/d})>\lambda\}\right)^{\frac{1}{q}}\\
  &\geq\frac{c_{\alpha,d}}{2\vartheta}\left(u\{|x|<x_{\vartheta}:\,\frac{c_{\alpha,d}}{\vartheta}(1-|x|^{\vartheta\alpha/d})>\frac{c_{\alpha,d}}{2\vartheta}\}\right)^{\frac{1}{q}}\\
  &=\frac{c_{\alpha,d}}{2\vartheta}u(B(0,x_{\vartheta}))^{\frac{1}{q}},
  \end{align*}
  where taking $x_{\vartheta}=(\frac{1}{2})^{d/\alpha\vartheta}$ in the last step. It now follows that for $0 <\vartheta< 1$,
  \begin{equation}\label{prove case 2_5}
  \|\mathcal{A}_{\alpha,1}^{\mathcal{S}}(u^{\frac{\alpha}{d}}f)\|_{L^{q,\infty}(u)}\gtrsim\frac{1}{\vartheta}\left(\frac{x_{\vartheta}}{\vartheta}\right)^{\frac{1}{q}}\simeq{\vartheta}^{-1-\frac{1}{q}}.
  \end{equation}
  Finally, combining $\eqref{A_1 constant}$, $\eqref{L^{p} value}$, $\eqref{prove case 2_5}$, and using that $\frac{1}{q}+\frac{\alpha}{d}=\frac{1}{p}$, we have that  $\eqref{prove case 2_4}$ holds for $\Phi(t)\geq ct^{1-\frac{\alpha}{d}}$, which gives the desired bound by the monotonicity of $\Phi$. $\hfill$ $\Box$

\end{document}